\documentclass[11pt,leqno]{article}
\usepackage{amsmath,amsthm,amsfonts,amssymb,amscd}
\input xypic
\usepackage[latin1]{inputenc}
\usepackage{amsmath}
\usepackage{amssymb}
\usepackage{amsthm}
\usepackage{graphicx}
\newtheorem{Corollary}{Corollary}
\newtheorem{Theorem}{Theorem}
\newtheorem{Proposition}{Proposition}
\newtheorem{Lemma}{Lemma}
\newtheorem{Claim}{Claim}
\theoremstyle{Definition}

\newtheorem{Example}{Example}
\theoremstyle{Remark}
\newtheorem{Remark}{Remark}

\def\bc{{\mathbb{C}}}

\def\fa{{\mathcal F}}

\def\Jac{\operatorname{{Jac}}}

\def\Area{\operatorname{{A}}}

\def\U{{\mathcal U}}
\def\al{{\alpha}}
\def\be{{\beta}}

\def\om{{\omega}}
\def\Om{{\Omega}}
\def\la{{\lambda}}

\def\ov{\overline}

\def\bq{{\Bbb Q}}

\def\re{{\Bbb R}}
\def\co{{\Bbb C}}

\def\po{{\partial}}

\def\vr{\varphi}

\def\sing{\operatorname{{Sing}}}
\def\sing{\operatorname{{sing}}}

\def\hot{\operatorname{{h.o.t.}}}
\def\Hol{\operatorname{{Hol}}}
\def\Diff{\operatorname{{Diff}}}

\def\S{\operatorname{{S}}}

\def\mfold{\operatorname{{M}}}

\def\sm{\setminus}

\bigskip

\title{Complex polynomial vector fields having a finitely curved orbit}
\author{Albet\~a C. Mafra}
\date{}

\begin{document}
\maketitle
\begin{abstract}
In this paper we address the following questions: (i) Let $C\subset
\bc^2$ be an orbit of a polynomial vector field which has finite
total Gaussian curvature. Is $C$ contained in an algebraic curve?
(ii) What can be said of a polynomial vector field which has a
finitely curved transcendent orbit? We give a positive answer to (i)
under some non-degeneracy conditions on the singularities of the
projective foliation induced by the vector field. For vector fields
with a slightly more general class of singularities we prove a
classification result that captures rational pull-backs of
Poincar\'e-Dulac normal forms.
\end{abstract}


 \section{Introduction and motivation}
\label{section:Intro} Consider a polynomial vector field
$X=P(x,y)\frac{\po}{\po x} + Q(x,y)\frac{\po }{\po y}$ on $\bc^2$.
The nonsingular orbits
 of $X$ are holomorphic curves in $\co^2$ and thus are oriented
 minimal surfaces in $\mathbb R^4$ with respect to the euclidian metric.
 They are Riemann surfaces, and their topology can be very complicated
(space of ends not denumerable, infinity genus and so on). In
general, for  a 2-dimensional oriented surface $L$ immersed in
$\re^4$, and parameterized locally by isothermal coordinates $z$,
that endow $L$ with a Riemann surface structure, we have a
naturally associated map $\Phi$, called the {\it holomorphic Gauss
map} of the immersion, from $L$ into the complex projective space
$\bc P^3$. If we consider $\bc  P^3$ with the homogeneous
coordinates $(z_1,z_2,z_3,z_4)$, the image of $\Phi$ is contained
in the quadric $\bq_2=\{z^2_1 +z^2_2+ z^2_3+z^2_4=0\}$
(\cite{Lawson}). In the case of a minimal immersion $\Phi$ is a
holomorphic map. The Gauss map is {\it algebraic} if $M$ is
conformally equivalent to a finitely punctured compact Riemann
surface $\ov M$ and $\Phi$ extends as a meromorphic map to $\ov
M$. A classical result due to Chern and Osserman  states this is
equivalent, for complete minimal immersions, to the finiteness of
the total curvature of the immersion. In particular, an algebraic
curve $C\subset \bc^2$ has algebraic Gauss map, nevertheless there
are holomorphic curves in $\bc^2$ with finite total curvature but
which are not algebraic curves, for instance orbits of suitable
vector fields (see Example~\ref{Example:Poincare-Dulac}). One of
the aims of this paper is to give conditions on the vector field
$X$ in order to assure that a nonsingular orbit with algebraic
Gauss map (i.e., with finite total curvature) is actually
algebraic. As we will see, this is related to the nature of the
singularities of the corresponding singular foliation on $\bc
P^2$. Let $\fa$ a holomorphic foliation with discrete singular set
$\sing(\fa)$ on a complex surface $M$.  A singularity $p\in
\sing(\fa)$ is called {\it irreducible} if there is an open
neighborhood $U$ of $p$ in $M$ where $\fa$ is induced by a
holomorphic vector field $Z$ of the form: $Z(x,y)=x\frac{\po}{\po
x} + [\la y +\hot ] \frac{\po}{\po y},\quad \la\notin \Bbb Q
_+\text{{\rm (non-degenerated)}}$, or
$Z(x,y)=x^{m+1}\frac{\po}{\po x} + [y(1+\la x^m) +
\hot]\frac{\po}{\po y} ,\quad m\geq 1 \text{{\rm (saddle-node)}}.$

Given a polynomial vector field $X$ with isolated singularities on
$\co^2$ we denote by $\fa(X)$ the corresponding holomorphic
foliation with singularities induced by $X$ on $\co P(2)$ (see
Example~\ref{Example:polynomial}). For the case where only
irreducible singularities  are allowed we have the following
statement:
\begin{Theorem}
\label{Theorem:combinatory} Let $X$ be a polynomial vector field
 on $\co ^2$ and let $L$ be an orbit of $X$ with finite
total curvature with respect to the metric induced by $\co^2$ on
$L$. Suppose that the singularities of $\fa (X)$ on $\bc P^2$ are
irreducible or that $X$ is without invariant lines on $\bc^2$. Then
$L$ is algebraic.
\end{Theorem}

In the second part of this work we classify polynomial vector
fields admitting a finitely curved transcendent orbit and whose
singularities of $\fa(X)$ are in the {\it Poincar\'e domain}
(\cite{Arnold}), i.e., $\fa(X)$ is given in a neighborhood of a
singular point by a germ of vector field of the form $Z=\lambda
x\frac{\partial}{\partial x} + \mu y\frac{\partial}{\partial
y}+..., \lambda/\mu \notin \mathbb R_-$. By Poincar\'e-Dulac
theorem such a singularity is {\it dicritical} if and only if it
is linearizable with $\lambda/\mu \in \mathbb Q_+$.   We prove:
\begin{Theorem}\label{Theorem:linearization}
 Let $X$ be a  polynomial vector field on $\co ^2$ such that the
singularities of  $\fa(X)$ are non-dicritical and in the
Poincar\'e domain. If $X$ has a finitely curved non-algebraic
orbit then $\fa(X)$ is given by a closed rational 1-form on $\bc
P^2$. Indeed, either $\fa(X)$ is a logarithmic foliation or there
is a rational map $f\colon \bc P^2 \dashrightarrow \bc P^2$ such
that $\fa(X)$ is the pull-back $f^*\fa(Y)$ where  $Y$ is a
Poincar\'e-Dulac normal form $Y(x,y)=(nx + cy^n)
\frac{\partial}{\partial x} + y \frac{\partial}{\partial y}$, for
some $n \in \mathbb N$ and some $c \in \bc \setminus \{0\}$. In
particular all orbits of $X$ have finite total curvature.
\end{Theorem}

\section{Preliminaries}
\label{section:Preli}

\subsection{Theorems of Chow and Remmert-Stein}
In this paragraph we introduce two extension theorems found in
\cite{[Gunning-Rossi]} that will be referred frequently in our
work; the theorem of Chow and the theorem of Remmert-stein.
\vglue.1in
 \noindent{\bf  Theorem of Remmert-Stein}. {\sl Let V be
subvariety of the polydisc D$\in \co^m$. Suppose $W$ is an
irreducible subvariety of $D\sm V$ of dimension $n$. If $n>dim(V)$,
then the closure $\ov{W}$ of $W$ in $D$ is an irreducible subvariety
of dimension $n$.} \vglue.1in

\noindent{\bf Theorem of Chow}. {\sl  If $X$ is closed and proper
analytic subset of $\co P(n)$, then $X$ is algebraic, i.e., $X$ is a
finite union of algebraic varieties. } \vglue.1in

\subsection{Holomorphic foliations with singularities}
Let $\mfold $ be a complex surface. By definition a codimension one
holomorphic foliation with singularities $\fa$ on $\mfold$ is given
by an open cover $\{U_j\}_{j\in J}$ of $\mfold$ and a collection of
1-forms $\{\om_j\}_{j\in J}$, where $\om_j$ is defined on $U_j$ and
such that if $U_j\cap U_k\neq \empty$, then
$$ \om_j|_{U_j\cap U_k}=g_{jk}\,.\,\om_k|_{U_j\cap U_k},$$
for some nonvanishing holomorphic function $g_{jk}$ defined in
$U_j\cap U_k$. The singular set of $\fa$ is the analytic subset
$\sing(\fa)$ of $\mfold$ defined by
$$\sing(\fa):=\{p\in U_j \,:\, \om_j(p)=0\}$$
It is well-known that $\sing(\fa)$ may be assumed of codimension
two. Thus $\sing(\fa)$ is a discrete subset of $\mfold$. The {\it
leaves} of $\fa$ are the leaves of the restriction $\fa |
_{\mfold\sm\sing(\fa)}.$ Given a singularity $p\in\sing(\fa)$, a
{\it separatrix} of $\fa$ at $p$ is a germ of analytic curve $\S$
at $p$ which is invariant and such that $p\in \S$. Such a curve
always exists (see section 1.4).
\begin{Example}[Polynomial vector fields]
\label{Example:polynomial} {\rm It is well-known that, from the
above definition,  a polynomial vector field
$X(x,y)=P(x,y){\po\over\po x}+Q(x,y){\po\over\po y}$ on $\bc^2$
induces a holomorphic foliation with singularities $\fa(X)$ on
$\bc P^2$. The foliation $\fa(X)$ is characterized by the fact
that if $L$ is a leaf of $\fa(X)$ which is not contained in the
line at the infinity then $L^*=L\cap \bc^2$ is a non-singular
orbit of $X$. Conversely, any foliation with singularities on $\bc
P^2$ is obtained in  this way. }
\end{Example}

\subsection{Resolution of singularities}
Given a holomorphic foliation with singularities  $\fa$ on a
complex surface $M$, a theorem of Seidenberg (\cite{Seidenberg})
gives a {\em resolution} of the singular points of $\fa$.
\vglue.1in \noindent{\bf Theorem of Seidenberg}. {\sl There is
finite sequence of blow-ups at the points of $\sing (\fa )$ such
that their composition gives a proper holomorphic map $\pi :\tilde
\mfold\rightarrow \mfold$ of a complex surface $\tilde \mfold$ and
a foliation $\fa^*=\pi^*(\fa )$ with isolated singularities such
that:
\begin{itemize}
\item[\rm (i)] $\pi^{-1}(\sing (\fa ))=\cup_{j=1}^{k=n} \Bbb P_j$
where each $\Bbb P_j$ is a projective line. $\pi^{-1}(\sing (\fa )$
is called the divisor of the resolution, and $\pi |_{\tilde \mfold
\sm \pi^{-1}(\sing (\fa )}$ is a biholomorphism.
\item[\rm(ii)] At any singularity $p\in\cup_{j=1}^{k=n} \Bbb P_j$
there is a local chart $(x,y)$ such that $x(p)=y(p)=0$ and $\fa^*$
is given by one of the 1-forms
$$ xdy- \la ydx+(\hot ),\quad \la\notin \Bbb Q _+\text{{\rm (non-degenerated)}},$$
$$x^{m+1}dy- y(1+\la x^m)dx+(\hot ),\quad m\geq 1 \text{{\rm
(saddle-node)}}.$$
\end{itemize}
}

\begin{Remark}{\rm
\begin{itemize}
\item[\rm (1)] Non-degenerated and saddle-node singularities are
usually called \textit{irreducible} singularities.
\item[\rm (2)] According to \cite{Camacho-Sad}, for each $p\in
\sing (\fa )$, $\fa$ admits at least one separatrix through $p$; if
the number of theses separatrices is finite, $p$ is called a {\em
non-dicritical} singularity and {\em dicritical} otherwise. The
singularity $p$ is  non-dicritical if and only if  all the
projective lines $\Bbb P_j$ belonging to $\pi ^{-1}(p)$ are
invariant by $\fa^*$.
\end{itemize}
}
\end{Remark}

\subsection{The Camacho-Sad Index Theorem}

Let $\fa $ be a holomorphic foliation on a complex surface $M$ as
above and $p\in \sing (\fa )$ an isolated singular point. Let $\S$
be a separatrix of $\fa$ at $p$. We can assume that $\S =\{q\quad
|\quad f(q)=0\},$ where $f$ is a holomorphic function defined in a
neighborhood $\U$ of $p$. We may assume that $f$ is reduced, that
is, $df\neq 0$ outside $p$. Then it is well-known that given a
holomorphic 1-form $\om$ defining $\fa $ in $\U$, with sing($\om$)=
there are holomorphic functions $g$ and $h$ in $\U$, and a
holomorphic 1-form $\eta$ in $\U$ such that  we have $ g\om
=kdf+f\eta$.
The \textit{Camacho-Sad index} is defined as
$$CS(\fa ,\S ,p):={1\over 2\pi i}\int_{\po\S}{\eta\over h}$$

\begin{Example} [Index of an irreducible singularity]
{\rm Let $\fa$ be a nondegenerate singularity given by the 1-form
$\om (x,y)=x(\la +yp(x,y))dy- y(\mu +xq(x,y))dx,$ with $\la .\mu
\neq 0$, and $p,q$ holomorphic functions. $\S_x=\{ (x,0)\quad
|\quad x\in\Bbb C\}$ and $\S_y=\{ (0,y)\quad |\quad y\in\Bbb C\}$
are separatrices for $\fa$. Then $CS(\fa ,\S_x ,0)= {\mu\over\la}$
and also $CS(\fa ,\S_y ,0)= {\la\over\mu}.$ \vglue.1in  Let now
$\fa$ be a saddle-node singularity given by the 1-form $\om
(x,y)=x^{n+1}dy- [y(1+\la x^n) +q(x,y)]dx,$ with $n\geq 1 $ and
$q$ holomorphic function. Let $\S$ be the strong separatrix for $
\om $
$$\S=\{ (0,y)|y\in\Bbb C\}.$$
 Then $CS(\fa ,\S ,0)= 0.$  If there is another separatrix $S'$ for
 $\fa$ then we can assume that $\om (x,y)=x^{n+1}dy- y[(1+\la x^n)
 +q(x,y)]dx$ and $S'=\{y=0\}$. For this we have
 $CS(\fa,S',0)=\lambda$.}
\end{Example}
The  Camacho-Sad Index theorem reads as follows
(\cite{Camacho-Sad}): \vglue.1in \noindent{\bf Index Theorem}. {\sl
Let $\S$ be a compact holomorphic curve in a complex surface
$\operatorname{M}$. Assume that $\S$ is invariant under a
holomorphic foliation $\fa$ with isolated singularities. Then the
number of singularities of $\fa $ on $\S$ is finite, and we have
$$\sum_{p\in \sing (\fa )\cap\S}CS(\fa ,\S_p ,p)=\S\cdot \S$$
where $\S_p$ is the germ of $\S$ at $p$  and  $\S\cdot\S$ is the
self intersection number  of the embedding $\S\hookrightarrow
\operatorname{M}$. } \vglue.1in

\subsection{The holomorphic Gauss map of a holomorphic vector field}
Let $X$ be a holomorphic vector field on $\bc^2$ and $L$ a
non-singular orbit of $X$. Since $L$ is a minimal immersion we can
define a holomorphic Gauss map $\Phi\colon L \to \mathbb C P(3)$
which takes values on the quadric ${\bq_2} : \sum\limits_{j=1}^4
z_j ^2=0$ (cf. \cite{Lawson}). This map is constructed as follows:
we define a {\it tangent  Gauss map\/} $\Phi\colon L \to G^o(2,4)$
from $L$ into the Grassmaniann of oriented two planes in $\re^4$
by $\Phi(p)=T_p(L)$.  Then using the canonical  analytical
isomorphism $G^o(2,4) \simeq {\bq_2 }$ we may consider the {\it
tangent Gauss map\/} $\Phi\colon L \to {\bq_2 } \subset \bc P^3$.
The map $\Phi$ is holomorphic and is called the {\it holomorphic
Gauss map} of the minimal immersion $L\subset \mathbb R^4$. Now we
exploit the fact that $L$ is a holomorphic curve tangent to $X$ in
$\bc^2$ to give another interpretation of the holomorphic Gauss
map. Given any $p\in L$ we have $T_p (L) =\bc . X(p)\subset
\bc^2$. Thus the map $\Phi$ takes values on the space of
directions on $\bc^2\setminus \{0\}$ which is naturally identified
with the projective line $\bc P^1$ the line at the infinity
$L_\infty=\bc P^2 \setminus \bc^2$. Therefore, we can consider the
holomorphic Gauss map of $L$ as a map $\Phi \colon L \to \mathbb
CP(1)$. The {\it spherical image} of $L$ is defined as
$\Phi(L)\subset \bc P^1$.

\subsection{The Integral Curvature Lemma}
\label{subsection:Integralcurvaturelemma}
 In this section we relate
the total curvature of a nonsingular orbit of a holomorphic vector
field in $\bc^2$ with the area of its image under the holomorphic
Gauss map. Let $\psi \colon L \to \mathbb R^n$ an oriented minimal
surface and let $\varphi \colon L\rightarrow \bc P^{n-1}$ be the
corresponding {\it holomorphic Gauss map} (see \cite{Lawson}). The
sphere $S^{2n-1}$ induces on $\bc P^{n-1}$ a metric with constant
holomorphic sectional curvature called the Fubini-Study metric,
which is given in homogeneous coordinates by
\[ ds^{2}=\frac{|Z\wedge dZ |^2 }{|Z|^4}.\]
For $\varphi (Z)=(\varphi_1 (Z),..., \varphi_n (Z))$ we can
calculate the area of $\Phi$ by evaluating the integral
\[ A(\Phi)=\int_L \Phi^* \omega \]
where $\omega$ is the area element induced by the Fubini-Study
metric restricted to the image of $\Phi$. The following result is
adapted from \cite{Santalo} Section 5   page 427:
\begin{Proposition}
\label{Proposition:integralcurvature} Let $\psi :L\longrightarrow
\mathbb{R}^n $ be a minimal immersion with holomorphic Gauss map
$\Phi$ and Gaussian curvature $K$, then
$$ \int_L K dA = -A(\Phi).$$
\end{Proposition}

Denote by $\sigma \colon \bc^2 \setminus \{0\} \to L_\infty$ the
canonical fibration where we identify $L_\infty \cong \mathbb C
P(1)$. Given a point $q\in L_\infty$ the fiber $\sigma^{-1}(q)$ is a
complex affine puncture punctured at one point. As a  corollary of
Proposition~\ref{Proposition:integralcurvature} we have:

\begin{Proposition}
\label{Proposition:curvaturearea} Let $L\subset \bc^2$ be a
nonsingular orbit of a holomorphic vector field $X$ on $\bc^2$.
Suppose that generically for $q \in L_\infty$ the intersection
number $\sharp\big(L \cap \sigma^{-1}(q)\big)$ is $\nu \in \mathbb N
\cup \{+\infty\}$. Then:
\begin{itemize}
\item[{\rm(i)}] If $\nu \in \mathbb N$ then the total curvature of
$L$ is finite equal to  $C(L)=-2\pi \nu$.
\item[{\rm(ii)}] If $\nu =+\infty$  then the total curvature of
$L$ is not finite.
\end{itemize}

\end{Proposition}

\begin{Remark}
{\rm One can deduce Propositions~\ref{Proposition:integralcurvature}
and ~\ref{Proposition:curvaturearea} from the {\it co-area formula}
that reads as follows: given two manifolds $M^m$ and $N^n$ with
$m\geq n$, a smooth map $\Phi :M\longrightarrow N$ and a measurable
function $ f:M\longrightarrow \mathbb{R}$ we have the following
$$ \int_M f(x)\|\Jac\Phi \|dx=\int_N \int_{\Phi^{-1}(y)}
f(z)d\mathcal{H}^{m-n}(z) dy  $$ where $\mathcal{H}^{m-n}$ is the
$(m-n)$-dimensional Hausdorff measure induced on $\Phi^{-1}(y)$ from
$M$, and $\|\Jac\Phi\|$ is the Radon derivative of the measure on
$N$ with respect to the image under $d\Phi$ of the $n$-dimensional
measure on $M$.}
\end{Remark}

\section{Examples}
\label{section:Examples}
 The total curvature of the orbits of a linear vector field is studied
 below.
\begin{Example}[Linear vector fields]
\label{Example:linear}{\rm  Let $X=x{\po\over \po x}+ \la y{\po
\over \po y}$ be a linear vector field on $\co^2$. A parametrization
for the orbit passing by $(1,y_0)\in \bc^2$ is given by $\psi
(z)=e^z{\po\over \po x}+ y_0e^{\la z}{\po \over \po y}.$ If
$\lambda\in \mathbb Q$ then clearly the orbits of $X$ are contained
in algebraic curves and therefore they have finite total curvature.
On the other hand for $\lambda$ not rational the map $\psi$ is
injective and computations inspired in
 \cite{Scardua} Section 5 show that the total curvature of the
orbit $L_{(1,y_0)}$ can be written as
$$C(\psi)=\int\int_{\mathbb R^2} {-2|\la |^2|y_0|^2(1+|\la
|^2-\la-\bar{\la}) e^{2(1+\al) u-2\be v}\over (e^{2u}+|\la
|^2|y_0|^2 e^{2\al u-2\be v})^2}du dv $$ Computations similar to
those in \cite{Scardua} {\rm(compare with Proposition 6, Section 5)}
actually show that $C(\psi)>-\infty$ if and only if $\la \in \Bbb
Q$}.
\end{Example}
 The following is a key example in our study.
\begin{Example}[Poincar\'e-Dulac normal form]
\label{Example:Poincare-Dulac}{\rm  Let $X=(nx+y^n){\po\over \po
x}+ y{\po \over \po y}$ be a Poincar\'e-Dulac normal form vector
field on $\co^2$ (\cite{Arnold}). For simplicity we will assume
that $n=1$.  The origin is a singularity exhibiting a unique
separatrix given by $\{y=0\}$. This is necessary for  the
existence of finitely curved orbits not contained in algebraic
curves as we will see in
Lemma~\ref{Lemma:localseparatrixaccumulated}. Let us show that
this is actually the case for the orbits of $X$. A parametrization
for the orbit passing by $(0,y_0)$ is given by $\psi
(z)=y_0ze^z{\po\over \po x}+ y_0e^z{\po \over \po y}.$
 Writing $\psi$ as map into $\re^4$, we get
$\psi (z)=({y_0ze^z+\bar{y_0}\bar{z}e^{\bar{z}}\over 2},
{y_0ze^z-\bar{y_0}\bar{z}e^{\bar{z}}\over
2i},{y_0e^z+\bar{y_0}e^{\bar{z}}\over 2},
{y_0e^z-\bar{y_0}e^{\bar{z}}\over 2i} ).$ The holomorphic Gauss map
of the immersion $\psi$ is given by $ \vr (z)={\po\psi\over \po z}=
({y(1+z)e^z\over 2}, {y(1+z)e^z\over 2i},{ye^z\over 2}, {ye^z\over
2i}).$
 Therefore denoting $z=u+iv$ and
$\la = \al +i\be$ we have the following expression for $F(z)
=|\vr|^2$
\begin{eqnarray*}
F(z) =  {1\over 4}|1+z|^2e^ze^{\bar{z}}+{1\over
4}|1+z|^2e^ze^{\bar{z}}+{1\over 4}|y|^2e^ze^{\bar {z}}+{1\over
4}|y|^2e^ze^{\bar {z}}  =  {|y|^2(|1+z|^2+1)e^{z+\bar{z}}\over 2}
\end{eqnarray*}
So $F(z,\bar{z})={|y|^2((1+z)(1+\bar{z})+1)e^{z+\bar{z}}\over 2}$,
\begin{eqnarray*}
{\po \over \po \bar{z}}\log F  = {[(1+z)(1+\bar{z})+1]+(1+z)\over
(1+z)(1+\bar{z})+1}+(1+z)+1 = {(1+z)(2+\bar{z})+1]+1\over
(1+z)(1+\bar{z})+1]+1}
\end{eqnarray*}
and
$${\po^2 \log(F) \over \po z\po\bar{z}}={1\over (|1+z|^2+1)^2}$$
Hence
\begin{eqnarray*}
K= -\frac{1}{F} {\po^2 \log(F) \over \po z\po\bar{z}}= {-2\over
|y|^2(|1+z|^2+1)e^{z+\bar{z}}}  {1\over (|1+z|^2+1)^2} ={-2\over
|y|^2(|1+z|^2+1)^3e^{z+\bar{z}}}
\end{eqnarray*}

In terms of the variables $u$ and $v$, $K$ can be written as
$$K={-2|\la |^2|y_0|^2(1+|\la |^2-\la-\bar{\la}) e^{2(1+\al)
u-2\be v}\over (e^{2u}+|\la |^2|y_0|^2 e^{2\al u-2\be v})^3}$$ For
the immersion $\psi (z)= (yze^z, ye^z)$ the induced metric is given
by
$$
ds^2= |y|^2e^{z+\bar{z}}(|1+z|^2+1)|dz|^2=
|y|^2e^{z+\bar{z}}(|1+z|^2+1)(|du|^2+|dv|^2).$$ The area element is
given by

$$dA =|y|^2e^{z+\bar{z}}(|1+z|^2+1)du\wedge dv$$
The total curvature can be written as
$$C(\psi) = \int\int _{\re^2}{-2e^{z+\bar{z}}(|y|^2e^{z+\bar{z}})(|1+z|^2+1)du
dv\over
|y|^2(|1+z|^2+1)^3}\\
= -2\int\int _{\re^2}{1\over (|1+z|^2+1)^2}du dv\\
= -2\pi $$
 Let us give a geometric interpretation of the above computation.
 The vector field  $X=(x+y) \frac{\po}{\po x} + y \frac{\po}{\po y}$
admits the first integral $f(x,y)=y e^{-\frac{x}{y}}$. Given $c\in
\bc\setminus \{0\}$ the orbit $L_c$ passing through the point
$(1,c)\in \bc^2$ is given by the level curve $f(x,y)=c
e^{-\frac{1}{c}}$. Given $\alpha \in \bc\setminus \{0\}$, the
intersection of the line $y=\alpha x$ with the leaf $L_c$ is given
by $f(x,\alpha x)=e^{-\frac{1}{c}}$ corresponds to $\alpha \, x \,
e^{-\frac{1}{\alpha}}= c. e^{-\frac{1}{c}}$ and therefore it is a
single point. Therefore, according to the notation of
Proposition~\ref{Proposition:curvaturearea} we have $\sharp(L_c \cap
\sigma^{-1}(q))=1$ for all $c\in \bc \setminus \{0\}$ and generic
$q\in L_\infty$. Applying this same proposition we conclude that the
area of the spherical image of $L_c$ is $2\pi$ which is the negative
of the total curvature of $L_c$. This is a general fact as follows
from Proposition~\ref{Proposition:curvaturearea}.

We point-out that the fact that $X$ has transcendent (non-algebraic)
orbits with  finite total curvature does not contradict
Theorem~\ref{Theorem:combinatory},
 indeed   this is due to the existence of  a  saddle-node
singularity in $L_\infty$. This saddle-node has strong manifold
contained in  the projective line $\ov\{y=0\}$ and central manifold
contained in the line $L_\infty$.

}
\end{Example}

\section{Finitely curved orbits of holomorphic vector fields}
\label{section:Finitelycurved}

 We shall now state two basic
properties of finitely curved orbits of holomorphic vector fields.
\begin{Lemma}[\cite{Scardua}]
\label{Lemma:flat} Let $L$ be an orbit of a holomorphic vector field
$X$ defined in $\co ^n$. If $L$ has finite total curvature as a real
surface defined in $\re ^{2n} $, then it is closed in $\co ^n
\setminus (\sing(X)\cup F(X))$, where $F(X)\subset \co ^n$ is the
union of the flat orbits of $X$.
\end{Lemma}
By a {\it flat orbit} we mean an orbit whose Gaussian curvature
vanishes identically. Such a orbit must be contained in a straight
complex line. Now we study the local behavior  of a finitely curved
orbit in a neighborhood of an irreducible singularity, obtaining the
following result:

\begin{Lemma}
\label{Lemma:localseparatrixaccumulated} Let $X$ be a holomorphic
vector field on $\bc^2$ and $p\in \sing(X)$ an irreducible
singularity. Let $L$ be an orbit of $X$ accumulating $p$ and if the
singularity $p$ is a saddle-node then its central manifold is also
accumulated by  $L$. If $L$ has finite total curvature then $L$ is
contained in the union of separatrices of $X$ through $p$.
\end{Lemma}
\begin{proof}   The central point is that, for an
irreducible  singularity, a non-separatrix leaf which accumulates on
the singularity and on the central manifold in the saddle-node case,
must accumulate on all the exceptional divisor after a blow-up at
the singular point.  By transverse uniformity, if $L$ accumulates a
point of a separatrix, it must accumulates all points in the
separatrix. The Gauss map for $L$ is given by $\Phi (p)=[X(p)],\quad
p\in L$ where $[(x,y)]$ denotes the straight line through the origin
passing through $(x,y)$ on $\co ^2$. It is then easy to see that the
map $\Phi$ is a holomorphic map from $L$ into $\co \Bbb P ^1$. As
usual we denote by $\fa(X)$ the corresponding holomorphic foliation
induced by $X$ on $\bc^2$. A blowing up at $p\in \co ^2$ will then
produce a foliation $\tilde \fa =\pi ^* (\fa(X))$ on $\tilde
{\bc}^2$. We can endow $\tilde{\co}^2$ with a Riemannian metric via
the pull-back $\pi ^* (ds^2 )$, where $ds^2$ is the canonical metric
on $\co ^2 $. The leaf $\tilde L = \pi ^{-1} (L)$ has the same
Riemannian metric behavior of $L$. So we get
$$\int_{\tilde L}\tilde K d \tilde A =\int_{L}K dA $$
In this case the projective line $\co P(1) $ is invariant with
respect to $\tilde \fa $, and since $\tilde L $ accumulates the
separatrices, $\tilde L $ must accumulate all points in $\co P(1) $.
It is also true that: outside  the separatrices, the Gauss map of
$\tilde L $ can be identified with the fiber map. Since $\tilde L $
accumulates $\co P(1)$ we have for $\tilde S=\sing(\tilde \fa)$ that
$ \# P^{-1}(p)\cap L^* =\infty $ for $p\in \co P(1)\sm \tilde S$ and
$P:\tilde{\co}^2 \longrightarrow \co P(1)$ being the fiber map.
Thus, $P$ has infinite area $A(P)=\infty$. Using the relation
$$\int_{L} K d  A=
\int_{\tilde L} \tilde K d \tilde A=-\Area(P)
$$ we get the desired conclusion. \end{proof}

 \begin{Corollary}
\label{Corollary:leafclosed} Let $X$ be a polynomial vector field
on $\bc^2$. Assume that every singularity $q\in \sing(\fa(X))$ is
irreducible. Then a finitely curved orbit $L$ of $X$ must be
closed outside $\sing(\fa(X))\cap \bc^2$.
 \end{Corollary}
\begin{proof}  If $L$ is not closed in $\bc^2\setminus
\sing(\fa(X))$ then $L$ accumulates some invariant straight line
$E\subset \bc^2$. By the Index Theorem this line $E$  contains
some singularity $q$ of $X$ in $\bc^2$ for which we have
$I(\fa(X),E,q) >0$. Therefore, if $q$ is a saddle-node then $E$ is
the central manifold of this saddle-node. Since the total
curvature of $L$ is finite  by
Lemma~\ref{Lemma:localseparatrixaccumulated} there is a
neighborhood $U\subset \bc^2$ of $q$ such that  $L\cap U$ is
contained in a local separatrix $\Gamma$ of $X$ through $q$. In
particular, $L$ cannot accumulate properly on the line $E$.
 \end{proof}

With the same proof of Lemma~\ref{Lemma:localseparatrixaccumulated}
we have:

\begin{Lemma}
\label{Lemma:Linfinity} Let $L$ be a nonsingular finitely curved
orbit of a polynomial vector field $X$ on $\bc^2$. Suppose that
the line at infinity is invariant by $\fa(X)$. If $L$ accumulates
on some point $p\in L_\infty$ then $p$ is a singularity.
\end{Lemma}

\section{Algebraicity of finitely curved orbits}
\label{section:algebraicity}
 In this section we prove
Theorem~\ref{Theorem:combinatory} and give a version of this result
for foliations with non-irreducible singularities but excluding the
saddle-node case (see Theorem~\ref{Theorem:finitecurvaturealgebraic}
below).

\begin{proof} [Proof of Theorem~\ref{Theorem:combinatory}] Let $X$ be a
polynomial vector field defined on $\co ^2$ whose corresponding
projective foliation  $\fa=\fa(X)$ has only irreducible
singularities and let $L$ be an orbit of $X$ with finite total
curvature. We study the behavior of $L$ in a neighborhood of
$L_\infty$. First we assume that the line $L_\infty$ is invariant.
We already know that $L$ is closed in $\bc^2 \setminus \sing(\fa)$
(Corollary~\ref{Corollary:leafclosed}). If $L$ accumulates some
regular point $p\in L_\infty$ then since $L_\infty$  is invariant
$L$ accumulates all the line $L_\infty$ and this is not possible
for $L$ has finite total curvature. Thus $L$ accumulates only at
singular points in $L_\infty$ and by Remmert-Stein theorem this
implies that $L$ has analytic closure of dimension one on $\bc
P^2$ and by Chow's Theorem $L$ is algebraic.
 Assume now that
$L_\infty$ is not invariant.  Let $p\in L_\infty$ be a regular point
accumulated by $L$. Then, since $L$ is closed in $\bc^2 \setminus
\sing(\fa)$ this implies that $L\cup \{p\}$ is analytic in a
neighborhood of $p$, indeed, $L\subset L_p$ the leaf through $p$. In
particular, $L$ accumulates only a discrete set of regular points in
$L_\infty\setminus \sing(\fa)$. Let now $p\in\sing(\fa)\cap
L_\infty$ be a singular point accumulated by $L$. Suppose that  $p$
is non-degenerate.  In this case $\fa$ exhibits two separatrices
through $p$ and either $L$ is contained in the union of separatrices
or $L$ accumulates both separatrices properly. The last possibility
implies that both separatrices are contained in (parallel) straight
lines in $\bc^2$ and $L$ accumulates both lines and is not closed in
$\bc^2\setminus \sing(\fa)$, absurd. Suppose now that  $p$ is a
saddle-node singularity. In this case either $L$ is contained in a
local separatrix of $\fa$ through $p$ or $L$ accumulates the strong
manifold of $\fa$ through $p$, this implies that $L$ is not closed
in $\bc^2 \setminus \sing(\fa)$, absurd.
 Since the set $\sing(\fa)\cap
L_\infty$ is also finite the set $L\cup \sing(\fa)$ is analytic and
therefore algebraic of dimension one on $\bc P^2$. This proves the
first part of the theorem. For the second part we assume that $X$
has no invariant line on $\bc^2$. This implies, as in the first
part, that $L$ is closed in $\bc^2\setminus \sing(\fa(X))$
(Lemma~\ref{Lemma:flat}). Now, essentially the same argumentation
above shows that $L$ must accumulate only on a finite set of points
in $L_\infty$. This shows, via the theorems of Remmert-Stein and
Chow, that $L$ is algebraic. This ends the proof of
Theorem~\ref{Theorem:combinatory}.
\end{proof}

\vglue.2in

According to \cite{C-LN-S} a singularity of a foliation in dimension
two is a {\it generalized curve} if its reduction process by
blow-ups exhibits no saddle-nodes. For this type of singularity
finitely curved leaves are always algebraic as follows:
\begin{Theorem}
\label{Theorem:finitecurvaturealgebraic} Let $X$ be a polynomial
vector field defined on $\co ^2$ and let $L$ be a finitely curved
 orbit of $X$.  Suppose that the singularities of
$\fa(X)$ on $\bc P^2$  are non-dicritical generalized curves then
$L$ is contained in an algebraic curve.
\end{Theorem}
The proof of Theorem~\ref{Theorem:finitecurvaturealgebraic} requires
the following lemma:

\begin{Lemma}
\label{Lemma:allseparatrices} Let $p\in\sing(\fa)$ be a
non-dicritical generalized curve. Then $\fa$ has at least two
separatrices through $p$. Moreover, if a leaf $L$ accumulates on $p$
and is not contained in the set of local separatrices of $\fa$
through $p$ then $L$ accumulates properly on all separatrices of
$\fa$ through $p$.
\end{Lemma}

The first part of Lemma~\ref{Lemma:allseparatrices} follows from
\cite{Mol}. The second part follows from the local study of
non-degenerate irreducible singularities and the Theorem of
Seidenberg (cf. Section 2.3).
\begin{proof} [Proof of
Theorem~\ref{Theorem:finitecurvaturealgebraic}] The proof has the
same structure of the proof of Theorem~\ref{Theorem:combinatory}.
First we prove that $L$ is closed in $\bc^2 \setminus \sing(\fa)$.
Indeed, if this is not the case, $L$ accumulates on a regular point
$q\in \bc^2 \setminus \sing(\fa)$ and therefore $L$ accumulates
properly on the leaf $L_q$ that must be an invariant straight line
in $\bc^2$. If $L_q$ contains some singularity $p\in \sing(\fa)\cap
\bc^2$ then by hypothesis $p$ is non-dicritical and it is a
generalized curve. By Lemma~\ref{Lemma:localseparatrixaccumulated}
we conclude that for some neighborhood $U$ of $p$ in $\bc^2$ the
intersection $L\cap U$ is contained in the set of local separatrices
of $\fa$ through $p$ and therefore the union $(L\cap U)\cup p$ is
analytic in $U$. This implies that $L$ cannot accumulate on the line
$L_q$ properly in $U$, contradiction. Therefore we must have
$L_q\cap \sing(\fa)\subset L_\infty$. In particular, either the
origin of the pencil $\sigma \colon \bc^2 \setminus \{0\} \to \bc
P^1\simeq L_\infty$ does not belong $L_q$ or it is not a singularity
of $\fa$. In both cases, since $L_q$ is invariant and accumulated by
$L$, the intersection number of $L$ with a generic fiber of $\sigma$
is infinity and therefore the total curvature of $L$ is infinite, a
contradiction. This proves that $L\cup \sing(\fa)$ is an analytic
subset of $\bc^2$. Now it remains to prove that $L$ accumulates only
on singular points in $L_\infty$.  If $L_\infty$ is invariant by
$\fa(X)$ this follows from Lemma~\ref{Lemma:Linfinity}. Assume now
that $L_\infty$ is not invariant. If there is a regular point $p\in
L_\infty$ which is accumulated by the orbit $L$ then the Flow Box
theorem shows that the leaf $L_p\not\subset L_\infty$ is properly
accumulated by $L$ in $\bc^2$, this is absurd because $L$ is closed
in $\bc^2 \setminus \sing(\fa(X))$. This proves
Theorem~\ref{Theorem:finitecurvaturealgebraic}.
\end{proof}

\section{Proof of Theorem~\ref{Theorem:linearization}}
\label{section:Prooflinearization}
 In this section we prove
Theorem~\ref{Theorem:linearization}. We shall use the following
proposition:
\begin{Proposition}
\label{Proposition:closedrational} Let $\fa$ be a holomorphic
foliation on $\bc P(2)$ and assume that:
\begin{itemize}
\item[{\rm (1)}] $\fa$ has an invariant {\rm(}irreducible{\rm)}
algebraic curve $\Lambda\subset \bc P(2)$. \item[{\rm (2)}] The
holonomy group of the leaf $\Lambda \setminus \sing(\fa)$ has  an
orbit accumulating only at the origin. \item[{\rm (3)}] The
singularities of $\fa$ on $\bc P(2)$ are non-dicritical and in the
Poincar\'e domain.
\end{itemize}
Then $\fa$ is given by a closed rational 1-form on $\bc P(2)$.
\end{Proposition}
\begin{proof}
Let $\fa$ be given in the affine space $\bc^2 \subset \bc P(2)$ by a
polynomial 1-form $\omega=Pdy - Qdx$ with isolated singularities.
The first step is to show the following
\begin{Claim}
The singularities of $\fa$ are either linearizable
{\rm(}non-resonant{\rm)} of the form $\lambda
x\frac{\partial}{\partial x} + \mu y\frac{\partial}{\partial y}$
with $\lambda /\mu \in \mathbb R \setminus \mathbb Q$ or
analytically conjugated to a Poincar\'e-Dulac normal form.
\end{Claim}
\begin{proof}
We take a singularity $p\in \Lambda \cap \sing(\fa)$. Since $p$ is
in the Poincar\'e domain we have two possibilities. Either $\fa$ is
analytically linearizable in a neighborhood of $p$ or $\fa$ is
analytically conjugated to a Poincar\'e-Dulac normal form in a
neighborhood of $p$. Moreover, since a Poincar\'e-Dulac normal form
exhibits only one separatrix, if there are two or more separatrices
then the singularity is analytically linearizable. On the other
hand, we are assuming that the singularities are non-dicritical,
therefore a linearizable singularity in the Poincar\'e domain cannot
be of resonant type, i.e., must be of the form $\lambda
x\frac{\partial}{\partial x} + \mu y\frac{\partial}{\partial y}$
with $\lambda /\mu \in \mathbb R \setminus \mathbb Q$.
\end{proof}

\begin{Claim}
The holonomy group $\Hol(\Lambda)$ is abelian. Moreover, either
$\sing(\fa)\cap \Lambda$ consists of only of linearizable
singularities or it consists only of Poincar\'e-Dulac type
singularities.
\end{Claim}
\begin{proof}
Indeed, the first remark is that by the Index theorem
$\sing(\fa)\cap \Lambda \ne \emptyset$. Given then a singularity
$p\in \Lambda\cap \sing(\fa)$, the local holonomy map $f\in
\Diff(\bc,0)$ of a local separatrix contained in $\Lambda$ gives
an element in the holonomy group $\Hol(\Lambda)$ of the leaf
$\Lambda \setminus \sing(\fa)$ which is either an analytically
linearizable non-periodic map of the form $z\mapsto e^{2\pi
\sqrt{-1}\lambda/\mu}$ with $\lambda / \mu \not \in \mathbb Q$, or
it is tangent to the identity and analytically conjugated to a map
of the form $z\mapsto \frac{z}{(z^n + 2\pi \sqrt{-1})^{1/n}}$. By
Nakai's Density theorem the holonomy group $\Hol(\Lambda)$ of the
leaf $\Lambda\setminus \sing(\fa)$ is solvable, maybe abelian. An
abelian subgroup of $\Diff(\bc,0)$ which contains a non-periodic
linearizable map is analytically linearizable. This implies that
the group $\Hol(\Lambda)$ is either abelian analytically
linearizable, abelian non-linearizable (without linearizable maps)
or solvable non-abelian and analytically conjugated to a subgroup
of a group $\mathbb H_k=\{(z\mapsto \frac{a z}{( bz^k +
1)^{1/k}})\}$ for some $k \in \mathbb N$ (\cite{Cerveau-Moussu}).
If the group $\Hol(\Lambda)$ contains a linearizable non-periodic
map and a map tangent to the identity conjugated  to an element of
$\mathbb H_k$ then its orbits are non-discrete off the origin
(\cite{Camacho-annals}) which is not the case. Therefore either
$\Hol(\Lambda)\setminus\{Id\}$ consists only of non-periodic
elements or it consists only of elements tangent to the identity.
In the first case, the group is necessarily abelian (because a
non-trivial commutator is tangent to the identity) and therefore
abelian linearizable. In this first case the singularities in
$\sing(\fa)\cap \Lambda$ are all linearizable. In the second case,
all singularities in $\sing(\fa)\cap \Lambda$ are of
Poincar\'e-Dulac type.  On the other hand in a solvable subgroup
of $\Diff(\bc,0)$ (maybe abelian) the subgroup of elements tangent
to the identity is  abelian. Thus we conclude that the group
$\Hol(\Lambda)$ is always abelian.
\end{proof}

Thus we have two possibilities for the holonomy group
$\Hol(\Lambda)$ and the singularities of $\fa$ in $\Lambda$.
Either $\Hol(\Lambda)$  is analytically linearizable and all
singularities of $\fa$ in $\Lambda$ are analytically linearizable
or $\Hol(\Lambda)$ is abelian non-linearizable and all
singularities of $\fa$ in $\Lambda$ are of Poincar\'e-Dulac type.
In the abelian linearizable case there is a closed meromorphic
1-form $\Omega$ with simple poles  defining $\fa$ in a
neighborhood of $\Lambda$ and by Levi's Extension theorem $\Omega$
extends to a closed rational 1-form of $\bc P(2)$ (see
\cite{Camacho-annals}, \cite{Scardua_Thesis}). This extension has
simple poles and can be written in a logarithmic form as
$\Omega\big|_{\bc^2}= \sum\limits_{j=1}^r
\alpha_j\frac{df_j}{f_j}$ for some irreducible polynomials $f_j$
and some complex numbers $\alpha_j\in \bc$. In this case the
foliation $\fa$ is a logarithmic or Darboux type foliation. Assume
now that the holonomy $\Hol(\Lambda)$ is  abelian
non-linearizable. In this case, we can once again construct a
closed meromorphic 1-form $\omega$ in a neighborhood of $\Lambda$
on $\bc P(2)$ (\cite{Scarduaparabolic}). This form is obtained as
follows: Fix a point $q\in \Lambda\setminus \sing(\fa)$ and choose
a transverse disc $\Sigma$ to $\fa$ at $q=\Sigma\cap \Lambda$ and
a local coordinate $z\in \Sigma$. The holonomy group
$\Hol(\Lambda)$ then corresponds to a subgroup
$\Hol(\fa,\Lambda,\Sigma)\subset \Diff(\bc,0)$ and the fact that
this group is abelian gives a germ of holomorphic vector field
$\xi(z)$ with a singularity at the origin $0\in \Sigma$ which is
invariant by the group $\Hol(\fa,\Lambda,\Sigma)$, i.e., $g_*
(\xi)=\xi, \, \forall g \in \Hol(\fa,\Lambda,\Sigma)$. Let
$\omega_0(z)$ be the 1-form dual to $\xi(z)$ in the sense that
$\omega_0(\xi)=1$. Then $\omega_0$ is also invariant by the
holonomy $\Hol(\fa,\Lambda,\Sigma)$. This invariance allows us to
extend $\omega_0$ by holonomy into a closed 1-form $\omega_1$ in a
neighborhood of $\Lambda\setminus \sing(\fa)$. It remains to show
that $\omega_1$ extends to a closed meromorphic 1-form $\omega$
which defines $\fa$ in a neighborhood of each singularity $p\in
\sing(\fa)\cap \Lambda$. This is done as in
\cite{Scarduaparabolic} as a consequence of the Poincar\'e-Dulac
normal form of these singularities (recall the fact that the
singularities of $\fa$ in $\Lambda$ are all Poincar\'e-Dulac of
same type). Again by Levi's Extension theorem the 1-form $\omega$
extends to a closed rational 1-form on $\bc P(2)$. This ends the
proof of the proposition.
\end{proof}

In the proof of Theorem~\ref{Theorem:linearization} the following
lemma will be useful:

\begin{Lemma}
\label{Lemma:PoincareDulac} Let $\fa$ be a holomorphic foliation
given by  a closed rational 1-form $\Om$ on  $\bc P(2)$. Assume
that the singularities of $\fa$ are
 non-dicritical in the Poincar\'e domain. Then either $\fa$ is a
 logarithmic foliation or it is a rational pull-back of a
 Poincar\'e-Dulac normal form.
\end{Lemma}
\begin{proof}
The main point is that, as we have observed above, a
non-dicritical singularity in the Poincar\'e domain is
linearizable if and only if it exhibits more than one local
separatrix. By the Integration Lemma (\cite{Scardua_Thesis}), in
an affine system of coordinates, $\Om$ can be written as
$$\Omega\big|_{\bc^2}=\sum\limits_{j=1}^r \lambda_j \frac{df_j}{f_j}
+ d\bigg(\frac{g}{\prod\limits_{j=1}^r f_j ^{n_j-1}}\bigg)$$ form
some irreducible polynomials $f_j$, some $n_j\in\mathbb N$ and  some
complex numbers $\lambda_j$. In the non-logarithmic case we have
some $n_j \geq 2$. Say for instance $n_1 \geq 1$. Suppose that
$\lambda_1 \ne 0$. We claim that $\lambda_j=0$ for all $j\geq 2$.
Indeed, if  for instance $\lambda_2 \ne 0$ then an intersection
point $q\in \ov\{f_1=f_2=0\}$ will exhibit two local separatrices
and therefore must linearizable, comparing this with the local form
of $\Omega$ at such a point we get a contradiction. The same
argumentation shows that $n_j=1$ for all $j\geq 2$. Thus we conclude
that $\Omega=\big|_{\bc^2}=\lambda_1 \frac{df_1}{f_1} +
d(\frac{g}{f_1 ^{n_1-1}})$ which is clearly a rational pull-back of
a Poincar\'e-Dulac normal form. Suppose now that $\lambda_1=0$. If
there is some $\lambda_j \ne 0$ we can assume that $\lambda_2 \ne
0$. In this case arguments as above show that $\lambda_j=0$ for all
$j\ne 2$ and we have $\Omega=\lambda_2 \frac{df_2}{f_2} + +
d\bigg(\frac{g}{\prod\limits_{j=1}^r f_j ^{n_j-1}}\bigg)$. Once
again we use the number of separatrices to show that we have a
contradiction with the fact that $n_1 \geq 2$. Therefore the only
possibility if that $\fa$ is a rational pull-back of a
Poincar\'e-Dulac normal form.
\end{proof}

\vglue.2in
\begin{proof}[Proof of Theorem~\ref{Theorem:linearization}]
 Let $L$ be a nonsingular transcendent orbit of $X$
with finite total curvature. By Lemma~\ref{Lemma:flat} we have two
possibilities:

 \noindent{\bf Case 1}. $L$ is closed in $\bc^2
\setminus \sing(\fa)$. In this case, since $L$ is not algebraic,
$L_\infty$ is invariant by $\fa$. Moreover,  given a small
transverse disc $\Sigma$ to $L_\infty$ at a point $q\in L_\infty
\setminus \sing(\fa)$, $L$ induces in $\Sigma$ an orbit which is
discrete outside the origin $q=\Sigma \cap L_\infty$. According to
Proposition~\ref{Proposition:closedrational}  $\fa(X)$ is given by
a closed rational 1-form $\Omega$ on $\bc P^2$. Using now
Lemma~\ref{Lemma:PoincareDulac} we conclude that $\fa$ is as
stated.

\noindent{\bf Case 2}.  $L$ is not closed in $\bc^2$ and  $L$
accumulates some invariant line $E\subset \bc^2$. The same
argumentation of the first case can be applied to $E$ in place of
$L_\infty$ to show that $\fa$ must be a rational pull-back of a
 Poincar\'e-Dulac normal form. This ends the proof
of Theorem~\ref{Theorem:linearization}.
\end{proof}

\bibliographystyle{amsalpha}

\begin{thebibliography}{31}
\frenchspacing
\bibitem{Arnold} V. Arnold, Chapitres Suppl\'ementaires de la
Th\'eorie des \'Equations Diff\'erentielles Ordinaires,  Mir 1980.

\bibitem{Camacho-Sad} C. Camacho, P. Sad,
{\it Invariants Varieties through singularities of holomorphic
vector fields}, Ann. of Math. 115 (1982), 579-595..

\bibitem{C-LN-S} C. Camacho, A. Lins Neto and P. Sad, {\it  Topological invariants and
equidesingularization for holomorphic vector fields}; J. of Diff.
Geometry, vol. 20, no.\ 1 (1984), 143--174.
\bibitem{Camacho-annals}
C. Camacho, A. Lins Neto and P. Sad, {\em  Foliations with algebraic
limit sets}; Ann. of Math. 136 (1992), 429--446.
\bibitem{Cerveau-Moussu} D. Cerveau et R. Moussu,
{\em  Groupes d'automorphismes de $({\Bbb C},0)$ et \'equations
diff\'erentielles $ydy+\cdots =0$}; Bull. Soc. Math. France, 116
(1988), 459--488.

\bibitem{[Gunning-Rossi]} R.C. Gunning, H. Rossi: Analytic
functions of several complex variables. \ Prentice Hall, N.J. 1965.

\bibitem{Lawson} B. Lawson,  {\it Lectures on minimal submanifolds}.
Vol. I. Second edition. Mathematics Lecture Series, 9. Publish or
Perish, Inc., Wilmington, Del., 1980

\bibitem{Mol} R. Mol, {\it Meromorphic first integrals: some extension results};
 Tohoku Math. J. (2) 54 (2002), no.
1, 85--104.

\bibitem{Nakai} I. Nakai, {\em Separatrices for nonsolvable dynamics on $C,0$}.
Ann. Inst. Fourier (Grenoble) 44 (1994), no. 2, 569--599.


\bibitem{Santalo} L. A. Santal\'{o}, {\it Integral geometry in Hermitian spaces.} Amer. J.
Math. 74, (1952). 423--434.



\bibitem{Scardua_Thesis} B. Sc\'ardua, {\it Transversely affine and
transversely projective holomorphic foliations having}; Ann. Sci.
Écolle Norm. Sup. vol. {\bf 30}, n 2, 1997.

\bibitem{Scardua} B. Sc\'ardua, {\it Complex vector fields having
orbits with bounded geometry}; Tohoku Math. J. (2) 54 (2002), no. 3,
367--392.


\bibitem{Scarduaparabolic} B. Sc\'ardua, \,{\it A Remark on Parabolic Projective Foliations},
Hokkaido Mathematical Journal, 28 (1999), 231--252.


\bibitem{Seidenberg} A. Seidenberg, {\it
Reduction of singularities of the differential equation $Ady=Bdx$};
Amer. J. of Math. 90 (1968), 248--269.






\end{thebibliography}

\begin{tabular}{ll}
Albet\~{a} Costa Mafra\\
 Instituto de  Matem\'atica - \\
 Universidade Federal do Rio de Janeiro\\
  Caixa Postal 68530\\
CEP. 21.945-970 Rio de Janeiro-RJ\\
 BRAZIL
 e-mail: albetan@im.ufrj.br
\end{tabular}

\end{document}